\documentclass[12pt]{article}
\usepackage{latexsym}
\usepackage{amsmath}
\usepackage{amssymb}

\usepackage{amsfonts}

\def\z{\zeta}

\newcommand{\C}{\mathbf{C}}
\newcommand{\Z}{\mathbf{Z}}
\newcommand{\eop}{\hfill$\square$}

\def\a{{\epsilon_p}}

\def\m1{{\epsilon_{\mu_1}}}
\def\m2{{\epsilon_{\mu_2}}}

\def\ea1{{\epsilon_{\mu_1^\ast}}}

\def\ea2{{\epsilon_{\mu_2^\ast}}}
\newtheorem{Con}{Conjecture}
\newtheorem{lemma}{Lemma}
\newtheorem{cor}[lemma]{Corollary}
\newtheorem{thm}[lemma]{Theorem}

\begin{document}

\centerline{{\Large \bf  Parametric Euler Sum Identities}}

\medskip
\centerline{David Borwein, Jonathan M.~Borwein, and David
M.~Bradley}

\medskip
\centerline{\today}

\bigskip {\bf Abstract.} We consider some parametrized classes of
multiple sums first studied by Euler.  Identities between
meromorphic functions of one or more variables in many cases
account for reduction formulae for these sums.

\bigskip

\section{\bf Introduction}
A somewhat unlikely-looking identity is
\begin{equation}
\label{eu-dual}
   \sum_{n=1}^{\infty}\frac1{n(n-x)}\sum_{m=1}^{n-1}\frac{1}{m-x}
    =\sum _{n=1}^{\infty}\frac {1}{n^2(n-x)},
\end{equation}
valid for all complex $x$ not a positive integer. For $x=0$,
\eqref{eu-dual} becomes $\zeta(2,1)=\zeta(3)$---the most central
Euler sum identity---as discussed below.

In this note  we begin with an \emph{ab initio} proof of
(\ref{eu-dual}), and  then explore  various consequences and
extensions. We conclude by studying other generating function
identities %parametric equations
of which
\begin{align}\label{gf2-3}
%\begin{split}
   \sum_{n=1}^\infty\frac{\displaystyle\frac{n}{{n}^{2}+{y}^{2}}
     +\sum_{m=1}^{n-1}\frac{2m}{m^2+y^2}}{n^2+4y^2}
   = \left(\coth(\pi y)+\coth(2\pi y)\right)
     \sum_{n=1}^{\infty}\frac{2\pi n y }{(n^2+y^2)(n^2+9y^2)}
%\end{split}
\end{align}
is a pretty example.

Though we do not belabour the point, all our work was assisted by
the use of  computer algebra systems---and some of it would have
been impossible, at least for us, without such tools. This is true
of both  our discoveries and of our proofs. The joys of such
symbolic and numeric computation $\ldots$ and much more are
discussed in detail in \cite{BBa} and \cite{BBaG}.

\section{A Parametric Euler sum}

We begin by defining various special functions which we shall
exploit in the sequel.  For ${\mathrm{Re}}(x)>0$, the gamma
function and its logarithmic derivative are defined by
\[
   \Gamma(x):=\int_0^\infty e^{-t}t^{x-1}\,dt,
   \qquad
   \Psi(x):=\frac {\Gamma'(x)}{\Gamma(x)}.
\]
Note that
\[
   \sum_{n=1}^\infty\frac{x}{n(n-x)}
   = -\Psi(1-x)-\gamma, \quad\mbox{where}\quad
   \gamma :=\lim_{N\to\infty}\bigg(\sum_{n=1}^N\frac1n-\log n\bigg)
\]
is Euler's constant.
%Hence, for $2\le N\in\Z$,
%\begin{equation}\label{diff-eq}
%   \sum_{n=1}^\infty\frac{x^N}{n^N (n-x)}
%   = -\zeta(N)x^{N-1} - \cdots-\zeta(2)\,x-\Psi(1-x)-\gamma.
%\end{equation}
%Similarly, with a little care integrating yields
%\begin{equation}\label{int-eq}
%   \sum_{n=1}^\infty\left\{{\frac{t}{n}}+\log \bigg( 1-{\frac {t}{n}}
%   \bigg)\right\}
%   = t\gamma-\log\Gamma(1-t).
%\end{equation}
%When $t=1/2$, Equation (\ref{int-eq}) evaluates to
%$(\gamma-\log\pi)/2$.
%\medskip

Recall that the classical\emph{ Riemann zeta  function} is defined
by
\[
   \zeta(s) :=  \sum_{n=1}^\infty \frac{1}{n^s},
   \qquad {\mathrm{Re}}(s)>1.
\]
Correspondingly,
\begin{equation}\label{doubleEuler}
   \zeta(s,t) := \sum_{n>m>0} \frac{1}{n^s m^t}
   =\sum_{n=1}^\infty \frac{1}{n^s}\sum_{m=1}^{n-1}\frac{1}{m^t}
%   =\sum_{n=1}^\infty \frac{H_{n-1}^{(t)}}{n^s},
\end{equation}
defines a double Euler sum.
%Here, $H_n^{(t)}:=1+1/2^t+\cdots+1/n^t$.
The two-place
function~\eqref{doubleEuler} was first introduced by Euler, who
noted the \emph{reflection formula}
\begin{equation}\label{refl}
   \zeta(s,t)+\zeta(t,s)=\zeta(s)\zeta(t)-\zeta(s+t),
   \qquad {\mathrm{Re}}(s)>1,\quad {\mathrm{Re}}(t)>1,
\end{equation}
and the reduction formula~\eqref{eu1} below.  An obvious extension
of~\eqref{doubleEuler} is
\[
   \zeta(s_1,s_2,\dots,s_N)
   :=\sum_{n_1>n_2>\cdots>n_N>0}\;
   \frac{1}{n_1^{s_1} n_2^{s_2}\cdots n_N^{s_N}},
\]
which defines an \emph{Euler sum} of \emph{depth} $N$ and
\emph{weight} $\sum_i s_i$.  Some authors reverse the order of the
variables.

Euler sums may be studied through a profusion of methods:
combinatorial, analytic and algebraic. The reader is referred to
\cite[Ch.\ 3]{BBaG} for a concise overview of Euler sums and their
applications.

We now prove identity~\eqref{eu-dual} directly.

\begin{thm} \label{th-main} If $x$ is any complex number not equal to a
positive integer, then
\[
   \sum_{n=1}^\infty\frac1{n(n-x)}\sum_{m=1}^{n-1}\frac1{m-x}
   =\sum_{n=1}^\infty\frac1{n^2(n-x)}.
\]
\end{thm}

\noindent{\bf Proof.} Fix $x\in\C\setminus\Z^{+}$. Let $S$ denote
the left hand side.   By partial fractions,
\begin{align*}
   S & %=\sum_{n=1}^\infty \sum_{m=1}^{n-1}\frac{1}{n(n-x)(m-x)}
   =\sum_{n=1}^\infty \sum_{m=1}^{n-1}\bigg(\frac{1}{n(n-m)(m-x)}
   -\frac{1}{n(n-m)(n-x)}\bigg)\\
   &=\sum_{m=1}^\infty\frac1{m-x}\sum_{n=m+1}^\infty\frac1{n(n-m)}
   -\sum_{n=1}^\infty \frac1{n(n-x)}\sum_{m=1}^{n-1}\frac1{n-m}\\
  &=\sum_{m=1}^\infty\frac1{m(m-x)}\sum_{n=m+1}^\infty\bigg(\frac1{n-m}
   -\frac1n\bigg)-\sum_{n=1}^\infty\frac1{n(n-x)}\sum_{m=1}^{n-1}\frac1m.
\end{align*}
Now for fixed $m\in\Z^{+}$,
\begin{align*}
   \sum_{n=m+1}^\infty \bigg(\frac1{n-m}-\frac1n\bigg)
   &= \lim_{N\to\infty}\sum_{n=m+1}^N \bigg(\frac1{n-m}-\frac1n\bigg)
   = \sum_{n=1}^m\frac1n-\lim_{N\to\infty}\sum_{n=1}^m\frac1{N-n+1}\\
   &= \sum_{n=1}^m \frac1n,
\end{align*}
since $m$ is fixed. Therefore, we have
\begin{align*}
   S &=\sum_{m=1}^\infty\frac1{m(m-x)}\sum_{n=1}^m\frac1n
   - \sum_{n=1}^\infty\frac1{n(n-x)}\sum_{m=1}^{n-1}\frac1m
   =\sum_{n=1}^\infty\frac1{n(n-x)}\bigg(\sum_{m=1}^n\frac1m-
   \sum_{m=1}^{n-1}\frac1m\bigg)\\
   &=\sum_{n=1}^\infty\frac1{n^2(n-x)}.
\end{align*}
\eop

\bigskip

For $x=\sqrt{-1}$, equation~\eqref{eu-dual} becomes the pair of
tangent sum evaluations
\begin{alignat*}{2}
   &\sum_{n=1}^\infty \frac1{n(n^2+1)}\sum_{m=1}^{n-1}\frac{mn-1}{m^2+1}
   &&=\sum_{n=1}^\infty \frac{1}{n(n^2+1)}
   = \gamma + {\mathrm{Re}}\,\Psi(1+\sqrt{-1})
\intertext{and}
   &\sum_{n=1}^\infty \frac1{n(n^2+1)}\sum_{m=1}^{n-1}\frac{n+m}{m^2+1}
   &&=\sum_{n=1}^\infty\frac{1}{n^2(n^2+1)}
   =\frac{\pi^2}6-{\rm Im}\,\Psi(1+\sqrt{-1})\\
   &&&=\frac{\pi^2}6-\frac{\pi\coth{\pi}-1}2.
\end{alignat*}
%where $\Psi$ denotes the logarithmic derivative of the Euler gamma
%function.

Setting $x=0$ in equation~\eqref{eu-dual} gives
$\zeta(2,1)=\zeta(3)$.  More generally,
differentiating~\eqref{eu-dual} $k$ times with respect to $x$
produces a corresponding formula for $\zeta(k+3)$. For example,
differentiating once and setting $x=0$ (equivalently, comparing
coefficients of $x$ on both sides) produces
\[
   \zeta(4)=\zeta(3,1)+\zeta(2,2).
\]
Since, by the \emph{reflection} formula (\ref{refl}), $\zeta(2,2)=
[\zeta^2(2)-\zeta(4)]/2$, we also evaluate
\[
  \zeta(3,1)=\frac 14\,\zeta(4)
  = \frac{\pi^4}{360}.
\]
This is the first case of a remarkable identity discovered by
Zagier, and first proved in~\cite{BBBL}.  See~\eqref{Zag} below.

\bigskip
Integrating equation~\eqref{eu-dual} produces the following
corollary.

\begin{cor}\label{log-cor} For all complex $x$ not equal to a
positive integer,
\begin{equation}\label{log}
   \sum_{n=1}^\infty\frac{1}{n^2}\,
   \log  \left( 1-\frac x n \right)
   = \sum_{n=1}^\infty\frac{1}{n}
   \sum_{m=1}^{n-1}\frac1{n-m}\log \left(\frac{1-x/m}{1-x/n}\right).
\end{equation}
\end{cor}
For $x=\sqrt{-1}$,  the imaginary part of~\eqref{log} leads to
\[
   \sum_{n=1}^{\infty}\frac1n
   \sum_{m=1}^{n-1}\frac1m\arctan\left(\frac{m}{n^2-mn+1}\right)
   =\sum_{n=1}^\infty{\frac{\arctan \left(1/n \right) }{n^2}}
   =\sum_{n=0}^\infty \frac{(-1)^n}{2n+1}\,\zeta(2n+3),
\]
where the last evaluation comes from writing
$\arctan(1/n)=\int_0^1 n(x^2+n^2)^{-1}\,dx$ and exchanging the
order of summation and integration.

\section{Euler's Reduction Formula}

Euler's  \emph{reduction formula} is
\begin{equation}\label{eu1}
   \zeta(s,1)
  =\frac{1}{2}s\zeta(s+1)
  -\frac12\sum_{k=1}^{s-2}\zeta(k+1)\zeta(s-k),
  \qquad 1<s\in\Z,
\end{equation}
which \emph{reduces} the double Euler sum $\zeta(s,1)$ to a sum of
products of classical Riemann $\zeta$-values.
%Euler also noted the first \emph{reflection formula}
%\begin{eqnarray}\label{refl}
%\zeta(a,b)+\zeta(b,a)=\zeta(a)\zeta(b)-\zeta(a+b),\end{eqnarray}
%certainly valid for $a,b>1$. This is an easy algebraic consequence
%of adding the double sums.
Another marvellous fact is the \emph{sum formula}
\begin{eqnarray}\label{sumf} \sum_{\substack{\Sigma a_i=s\\ a_i \ge
0}} \zeta\left(a_1+2,a_2+1,\cdots, a_r+1\right) =\zeta(r+s+1),
\end{eqnarray}
valid for all integers $s\ge 0$, $r\ge 1$. It is easy to show
that~\eqref{sumf} is equivalent to
\begin{equation}\label{sumg}
    \sum_{k_1>k_2>\cdots >k_r>0}\;\frac1{k_1}\prod_{j=1}^r\frac1{k_j-x}
     =\sum_{n=1}^\infty\frac{1}{n^r(n-x)},
     \qquad r\in\Z^{+}.
\end{equation}
The first three non-trivial cases of~\eqref{sumf} are
$\zeta(3)=\zeta(2,1)$, $\zeta(4)=\zeta(3,1)+\zeta(2,2)$ and
$\zeta(2,1,1)=\zeta(4)$.

The ordinary generating function of the sequence $\{\zeta(s) :
1<s\in\Z\}$ is
\[
   Z(x) := \sum_{s=2}^\infty \zeta(s) \,x^{s-1}
   = \sum_{n=1}^\infty \frac{x}{n(n-x)}
   = -\Psi(1-x)-\gamma.
\]
One can check that Euler's reduction~\eqref{eu1} is equivalent to
the ordinary generating function identity
\begin{equation}\label{zetas1GF}
   \sum_{n=1}^\infty \frac1{n(n-x)}\sum_{k=1}^{n-1}\frac1k
   =\sum_{n=1}^\infty\frac1{n(n-x)^2}
   -\frac{x}{2}\left\{\sum_{n=1}^\infty\frac1{n^2(n-x)^2}
   +\bigg(\sum_{n=1}^\infty\frac{1}{n(n-x)}\bigg)^2\right\}.
\end{equation}
This relies on observing that the right hand side of~\eqref{eu1}
involves the square and the derivative of $Z(x)$. In
turn,~\eqref{zetas1GF} is equivalent to~\eqref{eu-dual}.
%\begin{eqnarray}\sum _{n=1}^{\infty }{\frac {1}{{n}^{2} \left( n-x
%\right) }}=\sum _{n =1}^{\infty }{\frac {\sum
%_{m=1}^{n-1}\frac{1}{ \left( m-x \right) }}{n \left( n-x \right)}.
%}\label{eu-dual}\end{eqnarray}
This equivalence may be demonstrated as follows: Let
\begin{align*}
   &S_1 :=\sum_{n=1}^\infty\frac1{n(n-x)}\sum_{k=1}^{n-1}\frac1k,
   \quad
   &S_2 :=\sum_{n=1}^\infty\frac1{n(n-x)^2},
   \quad
   &S_3 :=\sum_{n=1}^\infty\frac1{n^2(n-x)^2},\\
   &S_4 :=\sum_{n=1}^\infty\frac1{n(n-x)},
   \quad
   &S_5 :=\sum_{n=1}^\infty\frac1{n^2(n-x)},
   \quad
   &S_6 :=\sum_{n=1}^\infty\frac1{n(n-x)}\sum_{k=1}^{n-1}\frac1{k-x}.
\end{align*}
It suffices to show that $S_1-S_2+ (S_3+S_4^2)x/2\equiv S_6-S_5.$
Observe that
\[
   S_1-S_6
   =-x\sum_{n=1}^\infty\frac{1}{n(n-x)}\sum_{k=1}^{n-1}\frac1{k(k-x)},
   \qquad
   S_5-S_2=-xS_3,
\]
and
\[
   S_4^2-S_3
   =2\sum_{n=1}^\infty\frac{1}{n(n-x)}\sum_{k=1}^{n-1}\frac1{k(k-x)}.
\]
Combining these yields the desired equivalence.

Correspondingly, equation~\eqref{eu-dual} is equivalent to
\[
   \zeta(s+3) =\sum_{\substack{a+b=s\\a,b\ge 0}} \zeta(2+a,1+b),
   \qquad 0\le s\in\Z,
\]
which is an inversion of Euler's reduction formula~\eqref{eu1},
and simultaneously recaptures the case $r=2$ of
equation~\eqref{sumf}.

Thus, Theorem~\ref{th-main} has established all of these results
directly, without using analytic methods, or sophisticated partial
fraction decompositions, as is usual. The first case, once more, is
$\zeta(2,1)=\zeta(3)$.
\bigskip

%\section{Some related identities}

\section{More Generating Functions}

Euler sums arise very naturally and it is perhaps a trick of time
that $\zeta(2n)$ is viewed as more natural than
\[
   \zeta(\{2\}_n) := \zeta(\underbrace{2,2,\dots,2}_n)
   =\sum_{k_1>k_2>\cdots>k_n>0}\;\prod_{j=1}^n
   \frac{1}{k_j^2}.
\]
Indeed, Euler's infinite product for the sine function is
precisely equivalent to
$$\frac{\sin(\pi x)}{\pi x} = \sum_{n=0}^\infty
(-1)^n\,\zeta(\{2\}_n)\,x^{2n},$$ from which we may deduce that
\[
   \zeta(\{2\}_n) = \frac{\pi^{2n}}{(2n+1)!},
   \qquad 0\le n\in\Z.
\]

Similarly, with $\omega$ a primitive cube root of unity,
\begin{align*}
   \sum_{n=0}^\infty \zeta(\{3\}_n)\,x^{3n}
   &=\prod_{n=1}^\infty\left(1+\frac{x^3}{k^3}\right)
    =\frac{1}{\Gamma(1+x)\Gamma(1+\omega x)\Gamma(1+\omega^2 x)}\\
   &= 1+\zeta(3)x^3+ \left(\frac12\,\zeta^2(3)
    -{\frac {8}{21}}\,\frac{\pi^6}{6!}\right)
   {x}^{6} +O\big(x^9\big),
\end{align*}
and allows one to show that $ \zeta(\{3\}_n)$ is always in the
ring generated by the numbers $\zeta(3k)$ ($k=1,2,\dots,n$).

By various methods, one can show that $$\zeta(\{3\}_n) =
\zeta(\{2,1\}_n)$$ for all $0\le n\in\Z$, while a proof of
\begin{equation}\label{open}
   \zeta(\{2,1\}_n) \stackrel{?}{=} 2^{3n}\,\zeta(\{\overline{2},1\}_n),
   \qquad 0\le n\in\Z
\end{equation}
remains elusive. The bar over the 2 on the right hand side
of~\eqref{open} signifies that in the summation, terms of the form
$1/k^2$ are to be multiplied by $(-1)^k$. Such sums are called
%\emph{multiple zeta values} or MZVs.
alternating Euler sums. Only the first case of~\eqref{open},
namely
\[
   \sum_{k=1}^\infty \frac{1}{k^2}\sum_{m=1}^{k-1} \frac 1m
   = 8\, \sum_{k=1}^\infty \frac{(-1)^k}{k^2}\sum_{m=1}^{k-1} \frac 1m
   \qquad(=\zeta(3))
\]
has a self-contained proof~\cite{BBaG,BBBL}. Indeed, the only
other proven case\footnote{This is an outcome of a complete set of
equations for  MZV's of depth four.} is
\begin{align*}
   \sum_{k=1}^\infty \frac{1}{k^2}\sum_{m=1}^{k-1}\frac 1m
   \sum_{p=1}^{m-1}\frac{1}{p^2}\sum_{q=1}^{p-1}\frac 1q
   &= 64\,\sum_{k=1}^\infty \frac{(-1)^k}{k^2}\sum_{m=1}^{k-1} \frac
   1m\sum_{p=1}^{m-1}\frac{(-1)^p}{p^2}\sum_{q=1}^{p-1}\frac 1q
   &(=\zeta(3,3)).
\end{align*}

There has been abundant evidence amassed to support~\eqref{open}
since it was first conjectured~\cite{BBB} in 1996. For example,
very recently Petr Lison\v ek checked the first $n\le 85$ cases to
1000 decimal places in about 41 hours with only the \emph{expected
roundoff error}. And he checked $n=163$ in ten hours.  If
true,~\eqref{open} would be the \emph{only} known identification
thus far of an infinite parametrized class of alternating Euler
sums with a corresponding class of non-alternating Euler sums. An
intriguing reformulation of~\eqref{open} is stated below.
\begin{Con} Define a sequence of polynomials $a_n=a_n(t)$ for positive
integers $n$ by $a_1=a_2=t^3$ and
\[
   n(n+1)^2a_{n+2} = n(2n+1)a_{n+1}+(n^3+(-1)^{n+1}t^3)a_n,
   \qquad n\ge 1.
\]
Then
\[
   \lim_{n\to\infty} a_n
   = t^3\prod_{n=1}^\infty \bigg(1+\frac{t^3}{8n^3}\bigg).
\]
\end{Con}

A more recondite generating function identity is
\[
   \sum_{n=0}^\infty\zeta(\{3,1\}_n)\,x^{4n}
   =\frac{\cosh(\pi x)-\cos(\pi x)}{\pi^2\,x^2},
\]
which is equivalent to Zagier's conjecture (subsequently proved):
\begin{equation}\label{Zag}
   \zeta(\{3,1\}_n)=\frac{2\pi^{4n}}{(4n+2)!},
   \qquad 0\le n\in\Z.
\end{equation}
The proof of~\eqref{Zag} (see~\cite[p.\
160]{BBaG},~\cite{BBBL},~\cite{BowBradRyoo}) devolves from a
remarkable factorization of the generating function in terms of
Gaussian hypergeometric functions:
\[
   \sum_{n=0}^\infty\zeta(\{3,1\}_n)\,x^{4n}
   = {}_2F_1(t,-t;1;1)\; {}_2F_1(it,-it;1;1),
\]
where $t=(1+i)x/2$.
%   = {\rm F}\left(x\frac{(1+i)}2,-x\frac{(1+i)}2;1;1\right)
%     {\rm F}\left(x\frac{(1-i)}2,-x\frac{(1-i)}2;1;1\right).$$

\section{Further extensions}
%For $a, b>0$ and integer, we define $$\zeta_a(s;x):=\sum_{n>0}
%\frac{1}{n^{a}(n-x)^s}, \qquad
% \zeta_b(s,t;x):=\sum_{n>0} \frac{1}{n^{b}(n-x)^s}\sum_{m=1}^{n-1} \frac{1}{(m-x)^t}.$$
% Moreover
%$$\zeta_{b+1}(s;x)=\zeta_b(s+1;x)- x\,\zeta_{b+1}(s+1;x),$$
%$$\zeta_{b+1}(s,t;x)=\zeta_b(s+1,t;x)-x\,\zeta_{b+1}(s+1,t;x).$$
%
%When $x=0$,  we reduce to $\zeta_a(s;0)=\zeta(a+s),$ and $
%\zeta_b(s,t;0)=\zeta(b+s,t).$
%
%
%Now equation~\eqref{eu-dual} may be written as
%$$\zeta_2(1;x)=\zeta_1(1,1;x).$$ Using $$\frac{d}{dx}
%\zeta_a(s;x)=s\,\zeta_a(s+1;x)$$ and $$\frac{d}{dx}
%\zeta_b(s,t;x)=s\,\zeta_b(s+1,t;x)+t\,\zeta_b(s,t+1;x)$$ we derive
%for $n=0,1,2,\ldots$ that
%\begin{eqnarray}\label{z2-id}
%\zeta_2(n+1;x)=\sum_{\stackrel{s+t=n}{s,t>0}} \zeta_1(1+s,1+t;x),
%\end{eqnarray}
%on applying Leibnitz' product rule to Equation (\ref{eu-dual}).
%
%
%\bigskip
%For sums of arbitrary depth $N$ we might define
% $$\zeta_b(\textbf{s};x):=\sum_{n_1>n_2\cdots >n_N>0} \frac{1}{n_1^{b}(n_1-x)^{s_1}\cdots
% (n_N-x)^{s_N}}.$$
%  In particular,  the
%generating function corresponding to
%$$\sum_{\stackrel{s_1+s_2+\cdots+s_N=w}{s_1>0,\cdots,s_N>0}} \zeta(1+s_1,s_2,\ldots,s_N)
%=\zeta(w+1)$$ is \begin{eqnarray}\label{zn-id}
%\zeta_2(w+1;x)=\sum_{\stackrel{\sum s_i=w}{s_i>0}}
%\zeta_1(1+s_1,\cdots,1+s_N;x).
%\end{eqnarray}
% For $N=2$, this is Equation (\ref{z2-id}).
%
%\bigskip
%
%\noindent \textbf{Problem 1.} \emph{Is this the appropriate
%generalization to arbitrary depth Euler sums? What about character
%sums ? }
%
%\bigskip
%
Relatedly and more centrally, Euler's reduction
formula~\eqref{eu1} has many extensions. For instance, with even
$a>0$  and odd $b>1$ with $a+b=2N+1$, one has~\cite{BBB}
\begin{align}
\begin{split}\label{zab}
   \zeta(a,b)
   &=\zeta(a)\zeta(b)+\frac12\bigg\{\binom{a+b}{a}-1\bigg\}
     \zeta(a+b)\\
   &-\sum_{r=1}^{N-1} \left\{\binom{2r}{a-1}+\binom{2r}{b-1}\right\}
     \zeta(2r+1)\zeta(a+b-1-2r),
\end{split}
\intertext{and hence by the reflection formula~\eqref{refl},}
\begin{split}\label{zba}
   \zeta(b,a)
   &= -\frac12\left\{1+\binom{a+b}{a}\right\}\zeta(a+b)\\
   &+\sum_{r=1}^{N-1}\left\{\binom{2r}{a-1}+\binom{2r}{b-1}\right\}
     \zeta(2r+1)\zeta(a+b-1-2r).
\end{split}
\end{align}
Although $\zeta(1)$ is undefined, Euler's reduction~\eqref{eu1}
implies that~\eqref{zab} also holds in the case $b=1$ if we
interpret $\zeta(1)=0$.  Thus,~\eqref{zab} is indeed an extension
of~\eqref{eu1}.

We next recast~\eqref{zab} using generating functions, and then
develop some further consequences and interesting special cases.
We initially keep $b=2t+1$ fixed and want to express
\[
   \sum_{N=t+1}^\infty \zeta(a,b)x^{2(N-t-1)}
   =\sum_{N=t+1}^\infty\sum_{n=1}^\infty
    \frac{x^{2(N-t-1)}}{n^{2(N-t)}}\sum_{m=1}^{n-1}\frac{1}{m^b}
   =\sum_{n=1}^\infty\frac{1}{n^2-x^2}\sum_{m=1}^{n-1}\frac{1}{m^b}.
\]
We suppose that $x^2$ is not a positive integer and set $D:=\dfrac
d{dx}.$  The requisite component generating functions are
\begin{align*}
   \sum_{N=t+1}^\infty\zeta(a)x^{2(N-t-1)}
   &=\sum_{N=t+1}^\infty\sum_{n=1}^\infty \frac{x^{2(N-t-1)}}{n^{2(N-t)}}
   =\sum_{n=1}^\infty\frac{1}{n^2-x^2},\\
   \sum_{N=t+1}^\infty\zeta(a+b)x^{2(N-t-1)}
   &=\sum_{N=t+1}^\infty\sum_{n=1}^\infty \frac{x^{2(N-t-1)}}{n^{2(N-t)+b}}
   =\sum_{n=1}^\infty\frac{1}{n^b(n^2-x^2)},\\
   \sum_{N=t+1}^\infty \zeta(a+b)\binom{a+b}b x^{2(N-t-1)}
   &=x^{-2}\sum_{N=t+1}^\infty\zeta(2N+1)\frac{D^b}{b!}x^{2N+1}\\
   &=\frac1{2x^2}\sum_{n=1}^\infty n\left(\frac1{(n-x)^{b+1}}+\frac{1}{(n+x)^{b+1}}\right)-\frac{\zeta(b)}{x^2},\\
   \sum_{N=t+1}^\infty x^{2(N-t-1)}\sum_{r=1}^{N-1}\binom{2r}{b-1} & \zeta(2r+1) \zeta(2N-2r)\\
   &=\bigg(\sum_{n=1}^\infty\frac1{n^2-x^2}\bigg)\sum_{n=1}^\infty\frac{D^{b-1}}{(b-1)!}\left(\frac{x^2}{n(n^2-x^2)}\right)\\
   &=\bigg(\sum_{n=1}^\infty\frac1{n^2-x^2}\bigg)\sum_{n=1}^\infty\frac12\left(\frac1{(n-x)^b}+\frac1{(n+x)^b}\right),
\end{align*}
and
\begin{align*}
   \sum_{N=t+1}^\infty x^{2(N-t-1)}&\sum_{r=1}^{N-1}\binom{2r}{a-1} \zeta( 2r+1) \zeta(2N- 2r)\\
   &=\sum_{m=1}^{t}\zeta(2m)x^{-1}\sum_{n=1}^\infty\frac{D^{b-2m}}{(b-2m)!}\left(\frac{x^2}{n(n^2-x^2)}\right)\\
   &=\sum_{m=1}^{t}\zeta(2m)x^{-1}\sum_{n=1}^\infty\frac12\left(\frac1{(n-x)^{b+1-2m}}-\frac1{(n+x)^{b+1-2m}}\right)\\
   &=\sum_{m=1}^{t}\zeta(b+1-2m)x^{-1}\sum_{n=1}^\infty\frac12\left(\frac1{(n-x)^{2m}}-\frac1{(n+x)^{2m}}\right).
\end{align*}
Combining these, with a fair amount of care, yields a generating
function for fixed $b$:
\begin{align*}
   \sum_{n=1}^\infty\frac1{n^2-x^2}\sum_{m=1}^{n-1}\frac{1}{m^b}&=\sum_{s=1}^\infty
   \zeta(2s,b)x^{2s-2}\\
   &=\zeta(2t+1)\sum_{k=1}^\infty\frac1{k^2-x^2} -\frac 12\,\sum _{n=1}^{\infty }{ \frac {1}{{n}^{2\,t+1}
      \left( {n}^{2}-{x}^{2} \right) }}\\
   &\quad+\frac 1{2}\, \sum _{n=1}^\infty \frac{n}{2x^2}\left\{ \frac{1}{
     \left( n-x \right) ^{2\,t+2}}+ \frac 1{  \left( n+x
     \right)^{2\,t+2}}-\frac 2 {n^{2\,t+2}} \right\}\\
   &\quad-\frac12\,\bigg(\sum_{n=1}^\infty\frac1{n^2-x^2}\bigg)\cdot
     \sum_{n=1}^\infty\left\{ \frac{1}{(n-x)^{2t+1}}+
     \frac1{(n+x)^{2t+1}}\right\}\\
   &\quad-\sum _{m=1}^{t}\zeta(2t+2-2m)
     \sum _{n=1}^{\infty }\frac 1{2x}\, \left\{ \frac1{(n-x)^{2m}}
      - \frac1{(n+x)^{2m}}\right\}.
      \end{align*}
Thus,
\begin{align*}
   \sum_{n=1}^\infty\frac1{n^2-x^2}\sum_{m=1}^{n-1}\frac{1}{m^b}&=
   \frac{\pi\cot\left(\pi x\right)}{4x}\sum_{n=1}^{\infty } \left\{
      \frac{1} {\left( n-x \right) ^{2\,t+1}}+ \frac{1}{ \left( n +x
      \right)^{2\,t+1}}-\frac2{n^{2t+1}}\right\}\\
   &-\frac12\,\sum_{n=1}^\infty\frac{n^{-2t-1}}{n^2-x^2}
   -\frac 1{2x}\,\sum _{m=1}^{t+1}\zeta(2t+2-2m)
     \sum_{n=1}^\infty\left\{\frac1{(n-x)^{2m}}
       - \frac1{(n+x)^{2m}}  \right\},
\end{align*}
using
\[
   \frac{\pi\cot(\pi x)}{2x}=\frac1{2x^2}-\sum_{k=1}^\infty\frac1{k^2-x^2}
   \quad \mbox{and}\quad \zeta(0)=-\frac12.
\]

Now we sum over the odd parameter, $b$, and obtain a two-variable
ordinary generating function expressible as:
\begin{align}\label{gf2} \nonumber
    \sum_{n=1}^\infty\frac1{n^2-x^2}\sum_{m=1}^{n-1} \frac{my}{m^2-y^2}
    &=\sum_{s>0,t\ge 0} \zeta(2s,2t+1)\,x^{2s-2}y^{2t+1} \\
    &=\frac\pi 2 x\cot(\pi x)
      \sum_{n=1}^\infty\frac{ny(4y^2-x^2)/(n^2-y^2)}
      {\left((n+x)^2-y^2\right)\left((n-x)^2-y^2\right)}
      \qquad   \nonumber \\
    &+\frac\pi 2 x\cot(\pi x)
      \sum_{n=1}^\infty\frac{ny}{\left((n+x)^2-y^2\right)
      \left((n-x)^2-y^2\right)}\qquad   \nonumber \\
    &+\frac{\pi}2 y\cot(\pi y)
       \sum_{n=1}^\infty\frac{2ny}
       {\left((n+x)^2-y^2\right)\left((n-x)^2-y^2 \right)}\nonumber\\
    &-\frac12\sum_{n=1}^\infty\frac{ny}
    {\left(n^2-y^2 \right)\left(n^2-x^2 \right)}.
\end{align}
Attractive specializations come with $x=\pm y$, $x=\pm 2iy$, $x=0$
and, after dividing by $y$, for $y=0.$ For example, at $(1/2,1/2)$
where the right hand side has a removable discontinuity, we deduce
that
\[
   \sum_{n=1}^\infty\frac{1}{4n^2-1}\sum_{m=1}^{n-1}\frac{m}{4m^2-1}
  =\frac{\pi^2}{64}-\frac{1}{16}.
\]

Note that the second, third  and fourth terms on the right hand
side of~\eqref{gf2} are expressible as
\begin{align*}
    \Sigma_2&:= -\frac{\pi}4 x\cot(\pi x)\,{\frac {\Psi(1+x+y) -\Psi(1+x-y) -\Psi(1-x+y) + \Psi(1-x-y) }{4xy}},\\
     \Sigma_3&:=\frac{\pi}2 y\cot(\pi y)\,{\frac {\Psi(1+x+y)-\Psi(1+x-y) -\Psi(1-x+y) + \Psi( 1-x-y) }{4xy}},\\
 \Sigma_4&:= {\frac {1-\left\{\Psi( 1+x) -\Psi(1-y) +\Psi(x)- \Psi( y) +\pi\cot(\pi x)\right\}
    y}{4(y^2-x^2)}},
\end{align*}
respectively while the first term has a corresponding evaluation.
It is
\begin{align*}
 \Sigma_1&:={y}^{2}\pi\cot\left(\pi\right)\,
 {\frac{\,\Psi( 1-x- y )+\Psi( 1+x-y)+\Psi(
  1+x+y) +\Psi( 1-x+y)
 }{8x( {x}^{2}-4\,{y}^{2}) }}\\
 &-y\pi\cot( \pi x) \,{\frac {\Psi( 1+x+y) -\Psi( 1+x-y) -\Psi( 1-x+y) +\Psi( 1-x-y)}{16({x}^{2}-4\,{y}^{2})}}\\
  & -{y}^{2}\pi\cot( \pi x)\,
{\frac { \Psi( y+1) +\Psi( 1-y) }{4x( {x}^{2}-4\,{y}^{2}) }}.
\end{align*}

We finish by expressing a form of this generating function concisely
as:
\begin{thm}\label{thm:gf2-2} For all $x$ and $y$ with squares not
equal to negative integers,
\begin{align*}%\label{gf2-2}
   \sum_{n=1}^\infty
   &\bigg(\frac{n}{n^2+y^2}+\sum_{m=1}^{n-1}\frac{2m}{m^2+y^2}\bigg)
   \frac1{n^2+x^2}\\
   &= (x^2-4y^2)\pi x \coth(\pi x)\sum_{n=1}^\infty
     \frac{n}{(n^2+y^2)\left\{(n^2-x^2+y^2)^2+(2nx)^2\right\}}\\
   &+( 2\pi y \coth(\pi y) + \pi x \coth(\pi x))\sum_{n=1}^\infty
     \frac{n}{(n^2-x^2+y^2)^2+(2nx)^2}.
\end{align*}
\end{thm}

\noindent {\bf Proof.} Replace $x$ by $ix$ and $y$ by $iy$
in~\eqref{gf2}, and then factor denominators and regroup the terms
as needed. \hfill{\bf QED}

\medskip

Letting $x$ and $y$ approach zero yields $\zeta(2,1)=\zeta(3)$
again. Setting $x=2y$ produces~\eqref{gf2-3}---the identity with
which we began.

 To conclude, we challenge the reader to explicitly obtain
the corresponding  two variable generating function for
$\z(2t+1,2s)$.
 \vfill

  \medskip
\noindent{\bf Acknowledgements.} J. Borwein's research was
supported by NSERC and the Canada
  Research Chair Programme, and D. Borwein's research was supported
  by NSERC.

\noindent Thanks are due to John Zucker and John Boersma for the
discussions which stimulated this work.

  \medskip
\noindent{\bf Key Words.} Euler sums, Zeta functions, Generating
functions, Multiple zeta values.

  \medskip
\noindent{\bf Classification Numbers.} \textbf{Primary} 33C99
\hfill  \textbf{Secondary} 11A99, 11M99

\vfill \noindent David Borwein \hfill{email: dborwein@uwo.ca}\\
Department of Mathematics\\University of Western Ontario\\
  London, Ontario N6A 5B7 Canada

\bigskip

\noindent Jonathan M. Borwein \hfill{email: jborwein@cs.dal.ca}\\
Faculty of Computer Science\\ Dalhousie University\\
  Halifax, NS  B3H~1W5 Canada

\bigskip
\noindent David M.~Bradley \hfill{email:
dbradley@math.umaine.edu}\\ Department of Mathematics \&
Statistics\\ University of Maine\\ 5752 Neville Hall\\ Orono,
Maine 04469-5752 U.S.A.

  \end{document}